\documentclass[12pt,a4paper,english,oneside,reqno]{amsart}
\usepackage[english]{babel}
\usepackage[T1]{fontenc}
\usepackage[ansinew]{inputenc}

\usepackage{amsmath}
\usepackage{amsfonts}
\usepackage{amssymb}
\usepackage{amsthm}
\usepackage{graphicx}
\usepackage{stmaryrd}
\usepackage{verbatim}
\usepackage{enumerate}
\usepackage{mathtools}
\usepackage{calrsfs}
\usepackage{hyperref}

\newcommand{\SR}[1]{}

\mathtoolsset{showonlyrefs=true}
\newtheorem{theorem}{Theorem}[section]
\theoremstyle{plain}

\newtheorem{corollary}[theorem]{Corollary}

\newtheorem{definition}[theorem]{Definition}

\newtheorem{lemma}[theorem]{Lemma}

\newtheorem{proposition}[theorem]{Proposition}

\theoremstyle{definition}
\newtheorem{remark}[theorem]{Remark}
\numberwithin{equation}{section}
\allowdisplaybreaks[1]
\begin{document}
\bibliographystyle{hplain}
\title[Radial multipliers on reduced free products]{Radial multipliers on reduced free products of operator algebras} 
\author{Uffe Haagerup}
\address{Department of Mathematical Sciences \\ University of Copenhagen \\ Universitetsparken 5 \\ 2100 Copenhagen Ø \\ Denmark}
\email{haagerup@math.ku.dk}
\author{Sören Möller}
\address{Department of Mathematics and Computer Science \\ University of Southern Denmark \\ Campusvej 55 \\ 5230 Odense M \\ Denmark}
\email{moeller@imada.sdu.dk}
\date{\today}

\begin{abstract}
Let $A_i$ be a family of unital C*-algebras, respectively, of von Neumann algebras and $\phi: \mathbb{N}_0 \to \mathbb{C}$. We show that if a Hankel matrix related to $\phi$ is trace-class, then there exists a unique completely bounded map $M_\phi$ on the reduced free product of the $A_i$, which acts as an radial multiplier. Hereby we generalize a result of Wysocza{\'n}ski for Herz-Schur multipliers on reduced group C*-algebras for free products of groups.
\end{abstract}

\maketitle



\section{Introduction}
Let $\mathcal{C}$ denote the set of functions $\phi$ on the non-negative integers $\mathbb{N}_0$ for which the matrix 
\begin{align*}
h = (\phi(i+j) - \phi(i+j+1))_{i,j \ge 0}
\end{align*}
is of trace class. Let $G = \ast_{i \in I} G_i$ be the free product of discrete groups $(G_i)_{i \in I}$. In \cite{Wysoczanski1995} , J. Wysocza{\'n}ski proved that if $\phi \in \mathcal{C}$ and $\tilde{\phi}: G \to \mathbb{C}$ is defined by $\tilde{\phi}(e)=\phi(0)$ and $\tilde{\phi}(g_1 \dots g_n) = \phi(n)$ for all $n > 0$ when $g_j \in G_{i_j}\setminus\{e\}$ and $i_1 \neq i_2 \neq \dots \neq i_n$, then $\tilde{\phi}$ is a Herz-Schur multiplier on $G$ and $\| \tilde{\phi }\|_{HS} \le \| \phi  \|_\mathcal{C}$, where $\| \cdot \|_\mathcal{C}$ is the norm on $\mathcal{C}$ defined in \eqref{eq_C_norm} below. In particular, there is a unique completely bounded map $M_\phi: C_r^*(G) \to C_r^*(G)$ such that $M_\phi(1)=\phi(0)1$ and 
\begin{align*}
M_\phi(\lambda(g_1 \dots g_n)) = \phi(n) \lambda(g_1 \dots g_n)
\end{align*}
when $g_j \in G_{i_j}\setminus\{e\}$ and $i_1 \neq i_2 \neq \dots \neq i_n$ as above, and $\|M\phi\|_{cb} \le \|\phi\|_\mathcal{C}$. Furthermore J. Wysocza{\'n}ski proved that $\|M_\phi\|_{cb} = \|\phi\|_\mathcal{C}$ in the cases when $|I|=\infty$ and $|G_i|=\infty$ for all $i \in I$. In the special case of $\phi_s(n)=s^n$ for $n \ge 0$ and $|s| < 1$ it follows that 
\begin{align*}
\|M_{\phi_s}\|_{cb} \le \|\phi_s\|_\mathcal{C} = \frac{|1-s|}{1-|s|}.
\end{align*}

In this paper we will show that every function $\phi$ from $\mathcal{C}$ gives rise to radial multipliers $M_\phi$ on reduced free products of $C^*$-algebras and reduced free products of von Neumann algebras (cf. Theorem \ref{thm_main_theorem}), satisfying $\|M_\phi\|_{cb} \le\|\phi\|_\mathcal{C}$. Radial multipliers of general reduced free products of $C^*$-algebras were first considered by {\`E}. Ricard and Q. Xu in \cite{RicardXu2006} and the weaker estimate $\|M_\phi\|_{cb} \le |\phi(0)| + \sum_{n=1}^\infty 4n|\phi(n)|$ can be obtained from \cite[Corollary 3.3]{RicardXu2006}.

The main result is proved in Section \ref{section_main_proof}. In Section \ref{section_examples} we discuss a related set of functions $\mathcal{C}'$ (cf. Definition \ref{def_class_Cprime}). It was used by T. Steenstrup, R. Szwarc and the first author in \cite{HaagerupSteenstrupSzwarc2010} to characterize radial multipliers on free groups $\mathbb{F}_n$ ($2 \le n \le \infty$). Moreover, C. Houdayer and {\`E}. Ricard used it in \cite{HoudayerRicard2011} to characterize multipliers on the free Araki-Woods factor $\Gamma(H_R,,U_t)''$ (cf. Section \ref{radial_subsection_RH_examples}).

In Section \ref{section_integral} we obtain an integral representation of functions in the class $\mathcal{C}$ which together with N. Ozawa's result in \cite{Ozawa2008}, shows that for every hyperbolic group $\Gamma$, and every $\phi \in \mathcal{C}$, the function
\begin{align}
\tilde{\phi}(x) = \phi(d(x,e))
\end{align}
is a completely bounded Fourier multiplier on $\Gamma$ (cf. Remark \ref{remark_phi_tilde_Fourier}).

\section{The main results}
\label{section_main_result}
We start by defining the class $\mathcal{C}$, crucial in what follows.
\begin{definition}
\label{def_class_C}
Let $\mathcal{C}$ denote the set of functions $\phi: \mathbb{N}_0 \to \mathbb{C}$ for which the Hankel matrix $h = (\phi(i+j)-\phi(i+j+1))_{i,j \ge 0}$ is of trace-class.
\end{definition}

If $\phi \in \mathcal{C}$, then $k = (\phi(i+j+1)-\phi(i+j+2))_{i,j \ge 0}$ is of trace-class, as well. Furthermore, we have
\begin{align}
\label{eq_sum_phi_n}
\sum_{n=0}^\infty \left| \phi(n)-\phi(n+1) \right| & \le \|h\|_1 + \|k\|_1 < \infty,
\end{align}
where $\|x\|_1 = Tr(|x|)$ is the trace-class norm for $x \in B(l^2(\mathbb{N}_0))$. This implies that $c = \lim_{n \to \infty} \phi(n)$ exists.
For $\phi \in \mathcal{C}$ set 
\begin{align}
\label{eq_C_norm}
\| \phi \|_\mathcal{C} = \|h\|_1 + \|k\|_1 + |c|.
\end{align}

The main result of this paper is the following generalization of Wysocza{\'n}ski's result:

\begin{theorem} \quad \\
\label{thm_main_theorem}
\begin{enumerate}
\item Let $\mathcal{A} = \ast_{i \in I} (\mathcal{A}_i,\omega_i)$ be the reduced free product of unital C*-algebras $(\mathcal{A}_i)_{i \in I}$ with respect to states $(\omega_i)_{i \in I}$ for which the GNS-representation $\pi_{\omega_i}$ is faithful, for all $i \in I$.

If $\phi \in \mathcal{C}$, then there is a unique linear completely bounded map $M_\phi: \mathcal{A} \to \mathcal{A}$ such that $M_\phi(1)= \phi(0) 1$ and $M_\phi(a_1 a_2 \dots a_n) = \phi(n) a_1 a_2 \dots a_n$ whenever $a_j \in \mathring{\mathcal{A}}_{i_j}=\operatorname{ker}(\omega_{i_j})$ and $i_1 \neq i_2 \neq \dots \neq i_n$.
Moreover $\| M_\phi \|_{cb} \le \|\phi\|_\mathcal{C}$.

\item Let $(\mathcal{M},\omega) = \bar{\ast}_{i \in I} (\mathcal{M}_i,\omega_i)$ be the w*-reduced free product of von Neumann algebras $(\mathcal{M}_i)_{i \in I}$ with respect to normal states $(\omega_i)_{i \in I}$ for which the GNS-representation $\pi_{\omega_i}$ is faithful, for all $i \in I$.

If $\phi \in \mathcal{C}$, then there is a unique linear completely bounded normal map $M_\phi: \mathcal{M} \to \mathcal{M}$ such that $M_\phi(1)= \phi(0) 1$ and $M_\phi(a_1 a_2 \dots a_n) = \phi(n) a_1 a_2 \dots a_n$ whenever $a_j \in \mathring{\mathcal{M}}_{i_j}=\operatorname{ker}(\omega_{i_j})$ and $i_1 \neq i_2 \neq \dots \neq i_n$.
Moreover $\| M_\phi \|_{cb} \le \|\phi\|_\mathcal{C}$.
\end{enumerate}
\end{theorem}

\begin{remark}
\label{radial_remark_equality}
By J. Wysocza{\'n}ski's result in \cite{Wysoczanski1995}, the norm estimates in Theorem \ref{thm_main_theorem} are best possible, as equality is attained if $|I| = \infty$ and $(A_i, \omega_i) = ( C^*_r(G_i), \tau_i)$ for a family $(G_i)_{i \in I}$ of infinite discrete groups, where $\tau_i$ is the canonical trace on $C^*_r(G_i)$ coming from the left regular representation. It would be interesting to know for which $(A_i, \omega_i)_{i \in I}$, the equality  $\|M_\phi\|_{cb} = \|\phi\|_\mathcal{C}$ holds for all $\phi \in \mathcal{C}$.
\end{remark}

We start by proving that the operator $M_\phi$ is unique, if it exists.

\begin{lemma}[Uniqueness in Theorem \ref{thm_main_theorem}]
The map $M_\phi$ is uni\-quely determined by the conditions in the hypothesis of Theorem \ref{thm_main_theorem}.
\end{lemma}

\begin{proof}
The algebra $\mathbb{C} 1 + \left( \sum_{i \in I} \mathring{\mathcal{A}}_i \right) + \left( \sum_{i_1 \neq i_2} \mathring{\mathcal{A}}_{i_1} \mathring{\mathcal{A}}_{i_2} \right) + \dots$ is norm dense in $\mathcal{A}$ and, respectively, $\mathbb{C} 1 + \left( \sum_{i \in I} \mathring{\mathcal{M}}_i \right) + \left( \sum_{i_1 \neq i_2} \mathring{\mathcal{M}}_{i_1} \mathring{\mathcal{M}}_{i_2} \right) + \dots$ is $\sigma$-weakly dense in $\mathcal{M}$.
As $M_\phi$ is bounded, it is then uniquely defined on all of $\mathcal{A}$, respectively, on all of $\mathcal{M}$.
\end{proof}

Now to prove Theorem \ref{thm_main_theorem} we start by showing that it is enough to prove the result for the special case of the algebras $M=B(H_i,\Omega_i)$ equipped with $\omega_i$, the vector state given by $\Omega_i$, as this will implie the result for general C*- and von Neumann-algebras.

\begin{proposition}
\label{prop_reduce_A_and_M_to_BH}
If Theorem \ref{thm_main_theorem} part (2) holds for $(\mathcal{M}_i,\omega_i) = (B(H_i), \omega_{\Omega_i})$ for Hilbert spaces $(H_i,\Omega_i)$ and associated vector states $\omega_i$ then Theorem \ref{thm_main_theorem} holds in general.
\end{proposition}

\begin{proof}
Assume Theorem \ref{thm_main_theorem} holds for $(B(H_i), \omega_{\Omega_i})$ for arbitrary $H_i$ and $\Omega_i$.
Now let $\mathcal{A} = \ast_{i \in I} (\mathcal{A}_i,\omega_i)$, respectively, $(\mathcal{M},\omega) = \bar{\ast}_{i \in I} (\mathcal{M}_i,\omega_i)$.
Let $(H_i,\Omega_i) = (H_{\omega_i},\xi_{\omega_i})$ be the Hilbert space and state coming from the GNS-representation of $\mathcal{A}_i$, respectively, $\mathcal{M}_i$, and let $(H,\Omega) = \ast_{i \in I} (H_i,\Omega_i)$ be their Hilbert space free product.

Now by \cite[Definition 1.5.1]{VoiculescuDykemaNica1992}, $\mathcal{A}_i$, respectively, $\mathcal{M}_i$ can be realized as subalgebras of $B(H)$ by the action defined as follows. If $a \in \mathcal{A}_i, \gamma_1 \otimes \dots \otimes \gamma_n \in H$ with $\gamma_j \in \mathring{H}_j := \Omega_j^\perp$ then
\begin{align}
a (\gamma_1 \otimes \dots \otimes \gamma_n) & = a(\Omega_i) \otimes \gamma_1 \otimes \dots \otimes \gamma_n
\end{align}
if $i \neq j$, and otherwise 
\begin{align}
a (\gamma_1 \otimes \dots \otimes \gamma_n) & = (a(\gamma_1)-\langle a(\gamma_1),\Omega_i \rangle \Omega_i) \otimes \gamma_2 \otimes \dots \otimes \gamma_n \\
& \qquad + \langle a(\gamma_1),\Omega_i \rangle \gamma_2 \otimes \dots \otimes \gamma_n.
\end{align}

Hence $M_\phi|_\mathcal{A}$ and $M_\phi|_\mathcal{M}$ can be obtained by restricting $M_\phi$ to the respective subalgebra of $B(H)$ and we then have $\| M_\phi|_\mathcal{A} \|_{cb} \le \| M_\phi \|_{cb} \le \|\phi\|_\mathcal{C}$, respectively, $\| M_\phi|_\mathcal{M} \|_{cb} \le \| M_\phi \|_{cb} \le \|\phi\|_\mathcal{C}$ which gives the desired general result.
\end{proof}

We will prove the special case considered in Proposition \ref{prop_reduce_A_and_M_to_BH} in the following sections.

\section{Preliminaries}
\label{section_preliminaries}

We start by introducing some notation.

Let $(H,\Omega)=\ast_{i \in I} (H_i, \Omega_i)$. Also denote $\mathring{H}_i = \Omega_i^\perp$, for $i \in I$. Then by the definition of the Hilbert space free product we have $H = \bigoplus_{n=0}^\infty H(n)$, where $H(n) := \bigoplus_{i_1 \neq \dots \neq i_n} \mathring{H}_{i_1} \otimes \dots \otimes \mathring{H}_{i_n}$, for $n > 0$, and $H(0)=\mathbb{C} \Omega$. We will denote the projection from $H$ to $H(n)$ by $P_n \in B(H)$, and let $Q_n := \sum_{k=n}^\infty P_k$.

Now choose orthonormal bases $\mathring{\Gamma}_i$ for $\mathring{H}_i$, then $\Gamma_i = \mathring{\Gamma}_i \cup \{ \Omega_i \}$ are bases for $H_i$.
Put $\Lambda(0) = \{ \Omega \}$ and $\Lambda(n) = \{ \gamma_1 \otimes \dots \otimes \gamma_n \colon \gamma_j \in \mathring{\Gamma}_{i_j}, i_1 \neq \dots \neq i_n \}$ for all $n \ge 1$.
Then $\Lambda(n)$ is an orthonormal basis for $H(n)$, for all $n \ge 0$ and $\Lambda = \bigcup_{n=0}^\infty \Lambda(n)$ is an orthonormal basis for $H$.
Note that $\Lambda(1) = \bigcup_{i \in I} \mathring{\Gamma}_i$ considered as a subset of $H$.

Now we can define the basic operators in $B(H)$.
Let $\gamma \in \Lambda(1)$. Let $L_\gamma, R_\gamma \in B(H)$ be the operators for which $L_\gamma \Omega = R_\gamma \Omega = \gamma$, and for $\chi = \chi_1 \otimes \dots \otimes \chi_n \in \Lambda(n)$ where $\chi_j \in \mathring{\Gamma}_{i_j}$ and $\gamma \in \mathring{\Gamma}_i$ we have
\begin{align}
L_\gamma (\chi) = \left\{ \begin{array}{cl}
\gamma \otimes \chi & \qquad \text{if $i \neq i_1$} \\
0 & \qquad \text{if $i = i_1$}
\end{array} \right.
\end{align}
respectively,
\begin{align}
R_\gamma (\chi) = \left\{ \begin{array}{cl}
\chi \otimes \gamma & \qquad \text{if $i \neq i_n$} \\
0 & \qquad \text{if $i = i_n$}.
\end{array} \right.
\end{align}

Note that $L_\gamma$ and $R_\gamma$ are well-defined partial isometries in $B(H)$. Moreover for all $\gamma \in \Lambda(1)$ and $n \ge 0$ we have $L_\gamma H(n) \subseteq H(n+1)$, respectively, $R_\gamma H(n) \subseteq H(n+1)$.
For $\gamma = \gamma_1 \otimes \dots \otimes \gamma_n \in \Lambda(n)$ denote $L_\gamma = L_{\gamma_1} L_{\gamma_2} \dots L_{\gamma_n}$, respectively, $R_\gamma = R_{\gamma_n} R_{\gamma_{n-1}} \dots R_{\gamma_1}$, where we set $L_\Omega = R_\Omega = 1$.

\begin{lemma}
\label{lem_dense_in_BHi}
Let $B(H_i)\mathring{} = \{ a \in B(H_i) \colon \langle a \Omega_i, \Omega_i \rangle = 0 \}$. Then the set $\operatorname{span}\left\{ \{ L_\gamma \colon \gamma \in \mathring{\Gamma}_i \} \cup \{ L_\gamma^* \colon \gamma \in \mathring{\Gamma}_i \} \cup \{ L_\gamma L_\delta^* \colon \gamma, \delta \in \mathring{\Gamma}_i \} \right\}$
is $\sigma$-weakly dense in $B(H_i)\mathring{}$ considered as a subset of $B(H)$.
\end{lemma}

\begin{proof}
Let $(e_{\gamma,\delta})_{\gamma,\delta \in \Gamma_i}$ be the matrix units of $B(H_i)$ corresponding to the basis $\Gamma_i$. Then $\operatorname{span}\{ e_{\gamma,\delta} \colon (\gamma,\delta) \neq (\Omega_i,\Omega_i) \}$
is $\sigma$-weakly dense in $B(H_i)\mathring{}$. Moreover, by the natural embedding of $B(H_i)$ in $B(H)$ one gets for $\gamma, \delta \in \mathring{\Gamma}_i$ that $L_\gamma = e_{\gamma,\Omega_i}$, $L^*_\gamma = e_{\Omega_i,\gamma}$, and hence $L_\gamma L^*_\delta = e_{\gamma,\delta}$, which proves the lemma.
\end{proof}

\begin{definition}
Let $a=(a_i)_{i \ge 0} \in l^\infty(\mathbb{N}_0)$. Denote by $D_a$ the operator which is defined by $D_a (\xi) = a_n \xi$ for $\xi \in \Lambda(n)$, $n \ge 1$, respectively, $D_a(\Omega) = a_0 \Omega$ and by linearity is extended to all of $H$.
\end{definition}

Note that $D_a = \sum_{n=0}^\infty a_n P_n$ and that $D_a$ is bounded with $\|D_a\| = \|a\|_\infty$.
Let $S$ denote the standard shift on $l^\infty(\mathbb{N}_0)$, i.e., for $a =(a_i)_{i \ge 0} \in l^\infty(\mathbb{N}_0)$, let $S(a_0,a_1,a_2, \dots ) = (0, a_0, a_1, a_2, \dots)$.
To ease notation, consider seperately the following two cases, which together contain all possible situations.

\begin{definition}
For $\xi \in \Lambda(k)$ and $\eta \in \Lambda(l)$, $k,l \ge 0$ we say that we are in
\begin{itemize}
\item \textbf{Case 1} if $\xi=\Omega$ or $\eta=\Omega$ or $k,l \ge 1$ and $\xi = \xi_1 \otimes \dots \otimes \xi_k$ and $\eta = \eta_1 \otimes \dots \otimes \eta_l$ where $\xi_k \in \mathring{\Gamma_i}, \eta_l \in \mathring{\Gamma_j}$ and $i \neq j$, $i,j \in I$,
\end{itemize}
respectively,
\begin{itemize}
\item \textbf{Case 2} if $k,l \ge 1$ and $\xi = \xi_1 \otimes \dots \otimes \xi_k$ and $\eta = \eta_1 \otimes \dots \otimes \eta_l$ where $\xi_k, \eta_l \in \mathring{\Gamma_i}$ for some $i \in I$.
\end{itemize}
\end{definition}

\section{Technical lemmas}
\label{section_technical_lemmas}

\begin{definition}
\label{def_Phi_abcd}
For $x,y \in l^2(\mathbb{N}_0)$ and $a \in B(H)$ set 
\begin{align}
\Phi^{(1)}_{x,y} (a) & := \sum_{n=0}^\infty D_{(S^*)^n x} a D^*_{(S^*)^n y} + \sum_{n=1}^\infty D_{S^n x} \rho^n(a) D^*_{S^n y},
\end{align}
respectively,
\begin{align}
\Phi^{(2)}_{x,y} (a) & := \sum_{n=0}^\infty D_{(S^*)^n x} a D^*_{(S^*)^n y} + \sum_{n=1}^\infty D_{S^n x} \rho^{n-1}(\epsilon(a)) D^*_{S^n y}
\end{align}
where $\rho(a) := \sum_{\gamma \in \Lambda(1)} R_\gamma a R^*_\gamma$ and $\epsilon(a) := \sum_{i \in I} q_i a q_i$ and $q_i$ is the projection onto $\operatorname{span}\{ \xi \in \Lambda(n): n \ge 1, \xi = \gamma_1 \otimes \dots \otimes \gamma_n, \gamma_n \in \mathring{\Gamma}_i \}$ for $i \in I$.
\end{definition}

\begin{lemma}
\label{lem_rhoLL_epsLL}
Let $k, l \ge 0$. Then for every $\xi \in \Lambda(k)$ and $\eta \in \Lambda(l)$ we have for all $n \ge 0$, $\rho^n( L_\xi L^*_\eta ) = L_\xi L^*_\eta Q_{l+n}$ and $\epsilon( L_\xi L^*_\eta ) = \rho( L_\xi L^*_\eta )$ in Case 1, and, respectively, $\epsilon( L_\xi L^*_\eta ) = L_\xi L^*_\eta$ in Case 2.
\end{lemma}

\begin{proof}
For the first statement observe that
\begin{align}
\rho^n( L_\xi L^*_\eta )
& = \sum_{\zeta \in \Lambda(n)} R_\zeta L_\xi L^*_\eta R^*_\zeta \\
& = L_\xi \left( \sum_{\zeta \in \Lambda(n)} R_\zeta R^*_\zeta \right) L^*_\eta \nonumber \\
& = L_\xi Q_{n} L^*_\eta \nonumber \\
& = L_\xi L^*_\eta Q_{l+n} \nonumber.
\end{align}

For the second statement, let $\chi \in \Lambda(m)$. If $m > l$ then
\begin{align}
\epsilon( L_\xi L^*_\eta )(\chi)
& = \sum_{i \in I} q_i L_\xi L^*_\eta q_i (\chi) \\
& = L_\xi \sum_{i \in I} q_i q_i L^*_\eta (\chi) \nonumber \\
& = L_\xi Q(1) L^*_\eta (\chi) \nonumber \\
& = L_\xi L^*_\eta Q(l+1) (\chi) \nonumber \\
& = L_\xi L^*_\eta (\chi) \nonumber
\end{align}

While if $m=l$, $\chi = \eta$ and $\eta_l \in \mathring{\Gamma}_{j}$ for some $j \in I$ we have
\begin{align}
\epsilon( L_\xi L^*_\eta )(\eta)
& = \sum_{i \in I} q_i L_\xi L^*_\eta q_i (\eta)
& = q_{j} L_\xi L^*_\eta q_{j} (\eta)
& = q_{j} L_\xi L^*_\eta (\eta)
& = q_{j} (\xi),
\end{align}
and this is equal to $0$ in Case 1 (i.e., $\xi_k \notin \mathring{\Gamma}_{j}$), respectively, equal to $\xi$ in Case 2 (i.e., $\xi_k \in \mathring{\Gamma}_{j}$).

Note that both sides vanish if $m=l$ and $\chi \neq \eta$, or $m < l$.
\end{proof}

To calculate the completely bounded norm of $\Phi^{(\cdot)}_{x,y}$ we use the following result from \cite{ChristensenSinclair1989}.

\begin{theorem}\cite[Theorem 1.3]{ChristensenSinclair1989}
\label{thm_Phi_cb}
If $ \| \sum_i u_i u_i^* \|, \| \sum_i v_i^* v_i \| < \infty$ for some $u_i, v_i \in B(H)$, then $\Phi(a) = \sum_i u_i a v_i $ defines a normal completely bounded operator on $B(H)$ and $\| \Phi \|_{cb} \le \| \sum_{i \in I} u_i u_i^* \| \|\sum_{i \in I} v_i^* v_i \|$.
\end{theorem}

Using this theorem we ontain the following cb-norm estimates.

\begin{lemma}
\label{lem_Phi1_norm}
For $x, y \in l^2(\mathbb{N}_0)$ we have $\| \Phi^{(1)}_{x,y} \|_{cb} \le \|x\|_2 \|y\|_2$, respectively, $\| \Phi^{(2)}_{x,y} \|_{cb} \le \|x\|_2 \|y\|_2$.
\end{lemma}

\begin{proof}
Let $\chi \in \Lambda(m)$. Then 
\begin{align}
\sum_{n=0}^\infty D_{(S^*)^n x} D_{(S^*)^n x}^* (\chi) 
& = \sum_{n=0}^\infty \overline{x(m+n)} D_{(S^*)^n x} (\chi) \\
& = \sum_{n=0}^\infty |x(m+n)|^2 (\chi) \nonumber \\
& = \left(\sum_{n=m}^\infty |x(n)|^2 \right) (\chi) \nonumber,
\end{align}
respectively,
\begin{align}
\sum_{n=1}^\infty \sum_{\zeta \in \Lambda(n)} (D_{S^n x} R_\zeta) (D_{S^n x} R_\zeta)^* (\chi) 
& = \sum_{n=1}^\infty \sum_{\zeta \in \Lambda(n)} D_{S^n x} R_\zeta R_\zeta^* D_{S^n x}^* (\chi) \\
& = \sum_{n=1}^{m} D_{S^n x} D_{S^n x}^* (\chi)
 = \sum_{n=0}^{m-1} |x(n)|^2 (\chi).
\end{align}
Here the second equality holds since $R_\zeta R_\zeta^*(\chi) = 0$ if $\zeta \neq \chi_{m-n+1} \dots \chi_m$ or if $n>m$. 

On the other hand,
\begin{align}
& \sum_{n=1}^\infty \sum_{\zeta \in \Lambda(n-1)} \sum_{i \in I} (D_{S^n c} R_\zeta q_i) (D_{S^n c} R_\zeta q_i)^* (\chi) \\
& = \sum_{n=1}^\infty \sum_{\zeta \in \Lambda(n-1)} \sum_{i \in I} D_{S^n c} R_\zeta q_i R_\zeta^* D_{S^n c}^* (\chi) \\
& = \sum_{n=1}^{m} \sum_{\zeta \in \Lambda(n-1)} D_{S^n c} R_\zeta R_\zeta^* D_{S^n c}^* (\chi) \\
& = \sum_{n=1}^{m} D_{S^n c} D_{S^n c}^* (\chi)
= \sum_{n=0}^{m-1} |c(n)|^2 (\chi) \nonumber.
\end{align}
The second equality holds since $\chi_{m-n+1} \in \mathring{\Gamma}_i$ (the rightmost element of $R_\zeta^*(\chi)$) for a unique $i \in I$ if $n-1 < m$, and there is no such $i \in I$ if $n-1 \ge m$. The third equality holds as $R_\zeta R_\zeta^*(\chi) = 0$ for $\zeta \neq \chi_{m-n+2} \dots \chi_m$. 
Hence
\begin{align*}
& \left( \sum_{n=0}^\infty D_{(S^*)^n x} D_{(S^*)^n x}^* + \sum_{n=1}^\infty \sum_{\zeta \in \Lambda(n)} (D_{S^n x} R_\zeta) (D_{S^n x} R_\zeta)^* \right) (\chi) \\
& = \left\|x \right\|_2^2 (\chi),
\end{align*}
respectively,
\begin{align*}
& \left( \sum_{n=1}^\infty D_{(S^*)^n x} D_{(S^*)^n x}^* + \sum_{n=1}^\infty \sum_{\zeta \in \Lambda(n-1)} \sum_{i \in I} (D_{S^n x} R_\zeta q_i) (D_{S^n x} R_\zeta q_i)^* \right) (\chi) \\
& = \left\|x \right\|_2^2 (\chi).
\end{align*}
Using these calculations for $x, y \in \in l^2(\mathbb{N}_0)$ and applying Theorem \ref{thm_Phi_cb} to get the desired result.
\end{proof}

\begin{lemma}
\label{lem_Phi1_LL_eq_sum_abLL}
Let $k,l \ge 0$. If $\xi \in \Lambda(k)$ and $\eta \in \Lambda(l)$ then
\begin{align}
\Phi^{(1)}_{x,y} (L_\xi L^*_\eta ) = \left( \sum_{t=0}^\infty x(k+t) \overline{y(l+t)} \right) L_\xi L^*_\eta
\end{align}
and, respectively, 
\begin{align}
\Phi^{(2)}_{x,y} (L_\xi L^*_\eta ) = \left\{ \begin{array}{ll}
\sum_{t=0}^\infty x(k+t) \overline{y(l+t)} L_\xi L^*_\eta & \qquad \text{in Case 1} \\
\sum_{t=0}^\infty x(k+t-1) \overline{y(l+t-1)} L_\xi L^*_\eta & \qquad \text{in Case 2}.
\end{array}
\right.
\end{align}
\end{lemma}

\begin{proof}
We prove this by showing that both sides act similarly on all simple tensors in $H$.

Indeed, let $m \ge 0$ and $\chi \in \Lambda(m)$ and let $n \ge 0$. If $\chi = \eta \otimes \zeta$, where $\zeta \in \Lambda(m-l)$ for some $l \ge 0$ we have for the common type of terms in $\Phi_{x,y}^{(1)}$ and $\Phi_{x,y}^{(2)}$ that
\begin{align}
\label{eq_DLLD_eq_xyLL}
D_{(S^*)^n x} L_\xi L^*_\eta D^*_{(S^*)^n y}(\chi) 
& = \overline{y(m+n)} D_{(S^*)^n x} L_\xi L^*_\eta  (\chi) \\
& = \overline{y(m+n)} D_{(S^*)^n x} (\xi \otimes \zeta) \nonumber \\
& = x(k+m-l+n) \overline{y(m+n)} L_\xi L^*_\eta (\chi) \nonumber. 
\end{align}
Otherwise, if there is no $\zeta$ such that $\chi = \eta \otimes \zeta$, then both sides vanish, wherein we have used the convention that $x(p)=0$ for $p<0$.

For the other type of terms in $\Phi_{x,y}^{(1)}$, we get
\begin{align}
\label{eq_DLrhoLD_eq_xyLL}
D_{S^n x} \rho^n(L_\xi L^*_\eta) D^*_{S^n y}(\chi) 
& = \overline{y(m-n)} D_{S^n x} \rho^n(L_\xi L^*_\eta)  (\chi) \\
& = \overline{y(m-n)} D_{S^n x} L_\xi L^*_\eta Q_{l+n} (\chi) \nonumber \\
& = \overline{y(m-n)} D_{S^n x} Q_{l+n} (\xi \otimes \zeta) \nonumber \\
& = x(k+m-l-n) \overline{y(m-n)} L_\xi L^*_\eta Q_{l+n}(\chi) \nonumber 
\end{align}
where in the second equality we use Lemma \ref{lem_rhoLL_epsLL} and the fact that both sides vanish if $n > m - l$.

We now estimate the other type of terms in $\Phi_{x,y}^{(2)}$. In Case 1 we similarly get
\begin{align}
\label{eq_DLrhoepsLD_eq_xyLL_case1}
& D_{S^n x} \rho^{n-1}(\epsilon(L_\xi L^*_\eta)) D^*_{S^n y}(\chi) 
= D_{S^n x} \rho^{n}(L_\xi L^*_\eta) D^*_{S^n y}(\chi) \\
& \qquad = x(k+m-l-n) \overline{y(m-n)} L_\xi L^*_\eta Q_{l+n}(\chi) \nonumber 
\end{align}
with both sides vanishing for $n > m-l$.

In Case 2 we get by Lemma \ref{lem_rhoLL_epsLL}
\begin{align}
\label{eq_DLrhoepsLD_eq_xyLL_case2}
& D_{S^n x} \rho^{n-1}(\epsilon(L_\xi L^*_\eta)) D^*_{S^n y}(\chi) 
= D_{S^n x} \rho^{n-1}(L_\xi L^*_\eta) D^*_{S^n y}(\chi) \\
& \qquad = x(k+m-l-n) \overline{y(m-n)} L_\xi L^*_\eta Q_{l+n-1}(\chi) \nonumber 
\end{align}
with both sides vanishing for $n > m-l+1$.

Combining \eqref{eq_DLLD_eq_xyLL} and \eqref{eq_DLrhoLD_eq_xyLL} we get
\begin{align}
\Phi^{(1)}_{x,y} (L_\xi L^*_\eta ) (\chi) 
& = \sum_{n=0}^\infty D_{(S^*)^n x} L_\xi L^*_\eta D^*_{(S^*)^n y}
+ \sum_{n=1}^\infty D_{S^n x} \rho^n(L_\xi L^*_\eta) D^*_{S^n y} \nonumber \\
& = \sum_{n=0}^\infty x(k+m-l+n) \overline{y(m+n)} L_\xi L^*_\eta (\chi) \nonumber \\
& \qquad + \sum_{n=1}^{m-l} x(k+m-l-n) \overline{y(m-n)} L_\xi L^*_\eta (\chi) \nonumber \\
& = \left( \sum_{n=l-m}^\infty x(k+m-l+n) \overline{y(m+n)} \right) L_\xi L^*_\eta (\chi) \nonumber \\
& = \left( \sum_{t=0}^\infty x(k+t) \overline{y(l+t))} \right) L_\xi L^*_\eta (\chi) \nonumber 
\end{align}
as desired.

Similarly in Case 1, combining \eqref{eq_DLLD_eq_xyLL} and \eqref{eq_DLrhoepsLD_eq_xyLL_case1} we get
\begin{align}
\Phi^{(2)}_{x,y} (L_\xi L^*_\eta ) (\chi) 
& = \left( \sum_{t=0}^\infty x(k+t) \overline{y(l+t)} \right) L_\xi L^*_\eta (\chi).
\end{align}

While in Case 2, combining \eqref{eq_DLLD_eq_xyLL} and \eqref{eq_DLrhoepsLD_eq_xyLL_case2} we get
\begin{align}
\Phi^{(2)}_{x,y} (L_\xi L^*_\eta ) (\chi)
& = \sum_{n=0}^\infty x(k+m-l+n) \overline{y(m+n)} L_\xi L^*_\eta (\chi) \nonumber \\
& \qquad + \sum_{n=1}^{m-l+1} x(k+m-l-n) \overline{y(m-n)} L_\xi L^*_\eta (\chi) \nonumber \\
& = \left( \sum_{t=0}^\infty x(k+t-1) \overline{y(l+t-1)} \right) L_\xi L^*_\eta (\chi) \nonumber.
\end{align}
This completes the proof.
\end{proof}

We now establish some technical results conserning maps $\phi \in \mathcal{C}$.

\begin{lemma}
\label{lem_properties_psi1_psi2}
Let $\phi \in \mathcal{C}$ and let $h, k$ and $c$ be as in Definition \ref{def_class_C}. Put $\psi_1(n) = \sum_{i=0}^\infty (\phi(n+2i)-\phi(n+2i+1))$ and $\psi_2(n) = \psi_1(n+1)$, for $n \ge 0$. Then $\phi(n) = \psi_1(n) + \psi_2(n) + c$ for $n \ge 0$ and the entries $h_{i,j}$ and $k_{i,j}$ of $h$ and $k$ are given by $h_{i,j} = \psi_1(i+j) - \psi_1(i+j+2)$, respectively, $k_{i,j} = \psi_2(i+j) - \psi_2(i+j+2)$, for $i,j \ge 0$.
\end{lemma}

\begin{proof}
By \eqref{eq_sum_phi_n} we have 
\begin{align}
| \lim_{n \to \infty} \psi_1(n) | & \le \lim_{n \to \infty} \sum_{i=0}^\infty |\phi(n+2i)-\phi(n+2i+1)| = 0.
\end{align}
A similar statement holds for $\psi_2$ and therefore $\lim_{n \to \infty} \psi_1(n) = 0$ and $\lim_{n \to \infty} \psi_2(n) = 0$. Next, let $n \ge 0$ be fixed. Then simple computations give $\psi_1(n) + \psi_2(n) = \phi(n) - c$, and $\psi_1(n) - \psi_1(n+2) = \phi(n)-\phi(n+1)$, respectively, $\psi_2(n)-\psi_2(n+2)=\phi(n+1)-\phi(n+2)$. Using these equations, we get the desired formulas for $h_{i,j}$, respectively $k_{i,j}$.
\end{proof}

\begin{remark}
\label{remark_hk}
Since $h$, $k$ are trace-class, it is well-known (cf. \cite[p. 13]{HaagerupSteenstrupSzwarc2010}) that there exist $x_i, y_i, z_i, w_i \in l^2(\mathbb{N}_0)$ such that $h = \sum_{i=1}^\infty x_i \odot y_i$ and $\sum \|x_i\|_2 \|y_i\|_2 = \|h\|_1$, respectively, $k = \sum_{i=1}^\infty z_i \odot w_i$ and $\sum \|z_i\|_2 \|w_i\|_2 = \|k\|_1.$ Here we use the notation $(u \odot v)(t) = \langle t,v \rangle u$, for $u,v,t \in l^2(\mathbb{N}_0)$.
\end{remark}

\begin{lemma}
\label{lem_psi_sum}
For $\psi_1$ and $\psi_2$ as in Lemma \ref{lem_properties_psi1_psi2}, and $x_i, y_i, z_i$, and $w_i$ as in Remark \ref{remark_hk} we have $\psi_1(k+l) = \sum_{i=1}^\infty \sum_{t=0}^\infty x_i(k+t) \overline{y_i(l+t)}$ and $\psi_2(k+l) = \sum_{i=1}^\infty \sum_{t=0}^\infty z_i(k+t) \overline{w_i(l+t)}$.
\end{lemma}

\begin{proof}
Let $k, l \ge 0$, then
\begin{align}
\psi_1(k+l) & = \sum_{t=0}^\infty \psi_1(k+l+2t) - \psi_1(k+l+2t+2) \\
& = \sum_{t=0}^\infty h_{k+t,l+t} \nonumber \\
& = \sum_{t=0}^\infty \sum_{i=1}^\infty x_i(k+t) \overline{y_i(l+t)} \nonumber \\
& = \sum_{i=1}^\infty \sum_{t=0}^\infty x_i(k+t) \overline{y_i(l+t)}  \nonumber 
\end{align}
where the sums are absolutly convergent, and we use Lemma \ref{lem_properties_psi1_psi2} for the first two equalities. A similar reasoning applies to $\psi_2$.
\end{proof}

\section{Proof of the main result}
\label{section_main_proof}

As shown in Section \ref{section_main_result}, it is enough to prove the following lemma in order to obtain the result of the main theorem.

\begin{proposition}
\label{prop_Mphi_exists_and_cb}
Let $(H,\Omega)=\ast_{i \in I} (H_i, \Omega_i)$ be the reduced free product of Hilbert spaces $(H_i)_{i \in I}$ with unit vector $\Omega_i$ and let $\omega_i(a) = \langle a \Omega_i,\Omega_i \rangle$ for $a \in B(H_i)$ where we realize $B(H_i)$ as subalgebras of $B(H)$ via the standard embedding from \cite[Definition 1.5.1]{VoiculescuDykemaNica1992}. Then for every $\phi \in \mathcal{C}$, there exists a linear completely bounded normal map $M_\phi$ on $B(H)$ such that $M_\phi(1) = \phi(0) 1$ and $M_\phi(a_1 a_2 \dots a_n) = \phi(n) a_1 a_2 \dots a_n$ whenever $n \ge 1$, $i_1, \dots i_n \in I$ with $i_1 \neq i_2 \neq \dots \neq i_n$ and $a_j \in B(H_{i_j})\mathring{} =  \text{ker}(\omega_{i_j})$. Moreover, $\|M_\phi\|_{cb} \le \|\phi\|_\mathcal{C}$.
\end{proposition}

The proof of Proposition \ref{prop_Mphi_exists_and_cb} will be divided into a series of lemmas. 

\begin{lemma}
\label{lem_equiv_T_mult_T_cases}
Let $T: B(H) \to B(H)$ be a bounded linear normal map, and let $\phi: \mathbb{N}_0 \to \mathbb{C}$. The following statements are equivalent.
\begin{enumerate}[(a)]
\item For all $n \ge 1$, $i_1, \dots i_n \in I$ with $i_1 \neq i_2 \neq \dots \neq i_n$ and $a_j \in B(H_{i_j})\mathring{} =  \text{ker}(\omega_{i_j})$, we have $T(1) = \phi(0) 1$ and $T(a_1 a_2 \dots a_n) = \phi(n) a_1 a_2 \dots a_n$.
\item For all $k,l \ge 0$ and $\xi \in \Lambda(k), \eta \in \Lambda(l)$ we have
\begin{align}
T( L_\xi L_\eta^* ) = \left\{ \begin{array}{ll}
\phi(k+l) L_\xi L_\eta^* & \qquad \text{in Case 1} \\
\phi(k+l-1) L_\xi L_\eta^* & \qquad \text{in Case 2}.
\end{array} \right.
\end{align}
\end{enumerate}
\end{lemma}

\begin{proof}
Assume (a) and let $k,l \ge 1$, $\xi \in \Lambda(k), \eta \in \Lambda(l)$. Now by the definition of $L_\xi$ we have $L_\xi L^*_\eta = L_{\xi_1} \dots L_{\xi_k} L_{\eta_l}^* \dots L_{\eta_1}^*$.

If we are in Case 1, there exist $i,j \in I, i \neq j$ such that $\xi_k \in \mathring{\Gamma}_i$ and $\eta_l \in \mathring{\Gamma}_j$. Hence all adjacent terms above are from different $B(H_i)\mathring{}$, hence $L_\xi L^*_\eta$ is of the form $a_1 \dots a_n$ in (a) with $n=k+l$.

On the other hand, if we are in Case 2, there exists $i \in I$ such that $\xi_k, \eta_l \in \mathring{\Gamma}_i$. In this case $L_{\xi_k} L_{\eta_l}^* \in B(H_i)\mathring{}$, hence $L_\xi L^*_\eta$ is of the form $a_1 \dots a_n$ in (a) with $n=k+l-1$.
Applying (a) we get the conclusion of (b) for $k,l \ge 1$.
If $k=0$ or $l=0$, e.g., $\xi=\Omega$ or $\eta=\Omega$ the result follows similarly by using $L_\Omega = 1$.

Assume (b). Using Kaplansky's density theorem \cite[Theorem 5.3.5]{KadisonRingrose1997} and the fact that the product is jointly $\sigma$-strong continuous on bounded sets, by Lemma \ref{lem_dense_in_BHi} it is enough to check that $T(1) = \phi(0) 1$ and $T(a_1 \dots a_n) = \phi(n) a_1 \dots a_n$ whenever $n \ge 1$ and $a_j \in \{ L_\gamma | \gamma \in \mathring{\Gamma}_{i_j} \} \cup \{ L_\gamma^* | \gamma \in \mathring{\Gamma}_{i_j} \} \cup \{ L_\gamma L_\delta^* | \gamma, \delta \in \mathring{\Gamma}_{i_j} \}$ where $i_1 \neq i_2 \neq \dots \neq i_n$.

It is easy to check that $L_\gamma^* L_\delta = 0$ when $\gamma, \delta \in \Lambda(1), \gamma \neq \delta$. In particular, $L_\gamma^* L_\delta =0$ when $\gamma \in \mathring{\Gamma}_{i_j}$ and $\delta \in \mathring{\Gamma}_{i_{j+1}}$, since $i_j \neq i_{j+1}$. Hence $a_1 \dots a_n = 0$, unless $a_1 \dots a_n = L_{\gamma_1} \dots L_{\gamma_k} L^*_{\delta_l} \dots L^*_{\delta_1}$ for some $\gamma_j \in \Gamma_{i_j}, \delta_s \in \Gamma_{r_s}$, $i_j \neq i_{j+1}, r_s \neq r_{s+1}$ and $i_1, \dots, i_k, r_1, \dots, r_l \in I$.

If we are in Case 1, we have $i_k \neq r_l$. Hence neighboring elements on the right hand side are from different $B(H_i)\mathring{}$ and thus $n=k+l$. If we are in Case 2, we have $i_k = r_l$. Hence $L_{\gamma_k} L^*_{\delta_l} \in B(H_{i_k})\mathring{}$, thus $n=k+l-1$. Now (b) gives the result for $k \ge 1$ or $l \ge 1$.
Moreover the $k=l=0$ case of (b) gives $T(1)=\phi(0) 1$.
\end{proof}

Next, we explicitly construct such a map $T$.

\begin{lemma}
\label{lem_T1_T2_psi1_psi2}
Let $\phi \in \mathcal{C}$. Define maps
\begin{align}
T_1 = \sum_{i=1}^\infty \Phi^{(1)}_{x_i,y_i} 
\quad \text{and} \quad
T_2 = \sum_{i=1}^\infty \Phi^{(2)}_{z_i,w_i}
\end{align}
where $\Phi^{(\cdot)}_{x,y}$ are as in Definition \ref{def_Phi_abcd}, $\psi_1, \psi_2$ as in Lemma \ref{lem_properties_psi1_psi2}, and $x_i, y_i, z_i, w_i$ as in Remark \ref{remark_hk}. Then for all $k,l \ge 0$, $\xi \in \Lambda(k), \eta \in \Lambda(l)$ we have $T_1( L_\xi L_\eta^* ) = \psi_1(k+l) L_\xi L_\eta^*$,
respectively,
\begin{align}
\label{eq_T2_psi2}
T_2(L_\xi L_\eta^*) = \left\{ \begin{array}{ll}
\psi_2(k+l) L_\xi L_\eta^* & \qquad \text{in Case 1} \\
\psi_2(k+l-2) L_\xi L_\eta^* & \qquad \text{in Case 2}.
\end{array} \right.
\end{align}
\end{lemma}

\begin{proof}
Let $k,l \ge 0$ and $\xi \in \Lambda(k), \eta \in \Lambda(l)$. Now by Lemma \ref{lem_Phi1_LL_eq_sum_abLL} and Lemma \ref{lem_psi_sum} we have
\begin{align}
T_1(L_\xi L^*_\eta) & = \sum_{i=1}^\infty \Phi^{(1)}_{x_i,y_i}(L_\xi L^*_\eta) \\
& = \left( \sum_{i=1}^\infty \sum_{t=0}^\infty x_i(k+t) \overline{y_i(l+t)} \right) (L_\xi L^*_\eta) \nonumber \\
& = \psi_1(k+l) L_\xi L^*_\eta \nonumber.
\end{align}

Furthermore, in Case 1 we have by Lemma \ref{lem_Phi1_LL_eq_sum_abLL} and Lemma \ref{lem_psi_sum}
\begin{align}
T_2(L_\xi L^*_\eta) & = \sum_{i=1}^\infty \Phi^{(2)}_{z_i,w_i}(L_\xi L^*_\eta) \\
& = \left( \sum_{i=1}^\infty \sum_{t=0}^\infty z_i(k+t) \overline{w_i(l+t)} \right) (L_\xi L^*_\eta) \nonumber \\
& = \psi_2(k+l) L_\xi L^*_\eta \nonumber,
\end{align}
respectively, in Case 2
\begin{align}
T_2(L_\xi L^*_\eta) & = \sum_{i=1}^\infty \Phi^{(2)}_{z_i,w_i}(L_\xi L^*_\eta) \\
& = \left( \sum_{i=1}^\infty \sum_{t=0}^\infty z_i((k-1)+t) \overline{w_i((l-1)+t)} \right) (L_\xi L^*_\eta) \nonumber \\
& = \psi_2(k+l-2) L_\xi L^*_\eta \nonumber.
\end{align}
This completes the proof.
\end{proof}

\begin{lemma}
\label{lem_T_LL_T_phi}
Define $T = T_1 + T_2 + c \operatorname{Id}$ where $\operatorname{Id}$ denotes the identity operator on $B(H)$. Then for $k,l \ge 0$ and $\xi \in \Lambda(k), \eta \in \Lambda(l)$ we have
\begin{align}
T( L_\xi L_\eta^* ) = \left\{ \begin{array}{ll}
\phi(k+l) L_\xi L_\eta^* & \qquad \text{in Case 1} \\
\phi(k+l-1) L_\xi L_\eta^* & \qquad \text{in Case 2.}
\end{array} \right.
\end{align}
\end{lemma}

Note that by Lemma \ref{lem_equiv_T_mult_T_cases} this implies that $T(1) = \phi(1) 1$ and that for $n \ge 1$,  $T_\phi(a_1 a_2 \dots a_n) = \phi(n) a_1 a_2 \dots a_n$.

\begin{proof}
Assume we are in Case 1, then
\begin{align}
T( L_\xi L_\eta^* ) & = T_1( L_\xi L_\eta^* ) + T_2( L_\xi L_\eta^* ) + c L_\xi L_\eta^*\\
& = \left(\psi_1(k+l) + \psi_2(k+l) + c \right) L_\xi L_\eta^* \nonumber \\
& = \phi(k+l) L_\xi L_\eta^* \nonumber.
\end{align}
Here we use the definition of $T$, then Lemma \ref{lem_T1_T2_psi1_psi2}, and lastly Lemma \ref{lem_properties_psi1_psi2}.
If we are in Case 2, we similarly get
\begin{align}
T( L_\xi L_\eta^* ) & = T_1( L_\xi L_\eta^* ) + T_2( L_\xi L_\eta^* ) + c L_\xi L_\eta^*\\
& = \left(\psi_1(k+l) + \psi_2(k+l-2) + c \right) L_\xi L_\eta^* \nonumber \\
& = \left(\psi_2(k+l-1) + \psi_1(k+l-1) + c \right) L_\xi L_\eta^* \nonumber \\
& = \phi(k+l-1) L_\xi L_\eta^* \nonumber.
\end{align}
Here we furthermore use $\psi_2(n)=\psi_1(n+1)$, for $n \ge 0$.
\end{proof}

By this result we have proven the existence of $M_\phi$ in Proposition \ref{prop_Mphi_exists_and_cb}, and it remains to calculate the cb-norm.

\begin{lemma}
\label{lem_T_cb_norm}
We have $\| T \|_{cb} \le \|\phi\|_{\mathcal{C}}$.
\end{lemma}

\begin{proof}
Let $x_i, y_i \in l^2(\mathbb{N}_0)$ then we have by Lemma \ref{lem_Phi1_norm} that $\| \Phi^{(1)}_{x_i,y_i} \|_{cb} \le \|x_i\|_2 \|y_i\|_2$.
Furthermore, since $T_1 = \sum_{i=1}^\infty \Phi^{(1)}_{x_i,y_i}$ we have by Remark \ref{remark_hk}
\begin{align*}
\| T_1 \|_{cb} \le \sum_{i=1}^\infty \| \Phi^{(1)}_{x_i,y_i} \|_{cb} \le \sum_{i=1}^\infty \|x_i\|_2 \|y_i\|_2 = \| h\|_1,
\end{align*}
respectively,
\begin{align*}
\| T_2 \|_{cb} \le \sum_{i=1}^\infty \| \Phi^{(2)}_{z_i,w_i} \|_{cb} \le \sum_{i=1}^\infty \|z_i\|_2 \|w_i\|_2 = \| k\|_1.
\end{align*}
Hence $\| T \|_{cb} \le \| T_1 \|_{cb} + \| T_2 \|_{cb} + \| c Id \|_{cb} \le \| h\|_1 + \| k\|_1 + |c| = \| \phi \|_\mathcal{C}$ as desired.
\end{proof}

Combining Lemmas \ref{lem_T_LL_T_phi} and \ref{lem_T_cb_norm} we obtain Porposition \ref{prop_Mphi_exists_and_cb}, and therefore an application of Proposition \ref{prop_reduce_A_and_M_to_BH} yields the conclusion of Theorem \ref{thm_main_theorem}.

\section{Examples}
\label{section_examples}

\subsection{\texorpdfstring{The case \boldmath$\phi_s(n) = s^n$}{A simple example}}
As a first example we will look at a simple $\phi$ where $\| \phi \|_\mathcal{C}$ can be calculated explicitly.

\begin{corollary}
\label{cor_phi_s}
Let $\mathbb{D} = \{ s \in \mathbb{C} | |s|<1 \}$ and $s \in \mathbb{D}$. Denote by $\phi_s$ the function $\phi_s(n) := s^n$. Then $\phi_s$ defines a radial multiplier $M_{\phi_s}$ on $\mathcal{A} = \ast_{i \in I} (\mathcal{A}_i,\omega_i)$, respectively, $(\mathcal{M},\omega) = \bar{\ast}_{i \in I} (\mathcal{M}_i,\omega_i)$ as in Theorem \ref{thm_main_theorem}. Moreover, $\| M_{\phi_s} \|_{cb} \le |1-s|/(1-|s|)$.
\end{corollary}

\begin{proof}
The conclusion follows from Theorem \ref{thm_main_theorem}, once we show that $\phi$ belongs to $\mathcal{C}$ and that $\| \phi \|_\mathcal{C} = \| h \|_1 + \| k \|_1 + |c| \le |1-s|/(1-|s|)$.

Observe first that $c = \lim_{n \to \infty} \phi(n) = \lim_{n \to \infty} s^n = 0$ as $|c|<1$. Furthermore, $\phi(i+j+1)-\phi(i+j+2) = s( \phi(i+j)-\phi(i+j+1) )$ so $k = s \cdot h$, hence $\| \phi \|_\mathcal{C} = (1+|s|) \| h \|_1$. Moreover $\phi(i+j)-\phi(i+j+1) = (1-s) s^{i+j}$ so $h = (1-s) m$, where $m$ is the matrix $m_{i,j} = s^{i+j}$. This gives $\| \phi \|_\mathcal{C} = (1+|s|) \| h \|_1 = (1+|s|)|1-s| \| m \|_1$. Now $m = a \odot \bar{a}$, where $a = (s^k)_{k \ge 0} \in l^2(\mathbb{N}_0)$, hence $\| m \|_1 = \|a\|_2^2 = 1/(1-|s|^2)$. Combining these calculations we get $\| \phi \|_\mathcal{C} = |1-s|/(1-|s|)$, which proves the corollary.
\end{proof}

\subsection{Wysocza{\'n}ski's theorem}
As a second example, we will show that Wysocza{\'n}ski's result, apart from determining when equality holds, is a special case of Theorem \ref{thm_main_theorem}.

\begin{theorem}[{\cite[Theorem 6.1]{Wysoczanski1995}}]
\label{thm_Wysoczanski2}
Let $G = \ast_{i \in I} G_i$ be the free product of a family of discrete groups, and let $g \in G$ be $g = g_1 g_2 \dots g_n$ where $g_j \in G_{i_j} \backslash \{e\}$, $j_1, \dots, j_n \in I$ and $j_1 \neq j_2 \neq \dots \neq j_n$. If $\phi \in \mathcal{C}$ then  $\tilde{\phi}(g) = \phi(n)$ is a Herz-Schur multiplier on $G$. Moreover $\| \tilde{\phi} \|_{HS} \le \| \phi \|_\mathcal{C}$.
\end{theorem}

\begin{proof}
Let $\phi \in \mathcal{C}$ and $g = g_1 \dots g_n \in G$ as above. Now by \cite[p. 301]{BozejkoFendler1984} and \cite{deCanniereHaagerup1985} we have $\| \tilde{\phi} \|_{HS} = \| \tilde{\phi} \|_{M_0A(G)} =  \| \tilde{M}_{\tilde{\phi}} \|_{cb}$ where $\tilde{M}_{\tilde{\phi}}$ is the operator $\tilde{M}_{\tilde{\phi}}(\lambda(g)) = \tilde{\phi}(g) \lambda(g)$ for $g \in G$ and $\lambda$ the left regular representation. By the definition of $\tilde{\phi}$ this is $\tilde{M}_{\tilde{\phi}}(\lambda(g)) = \phi(n) \lambda(g)$ and by the definition of $L(G)$ we have $\lambda(g) = \lambda(g_1) \lambda(g_2) \dots \lambda(g_n)$. Hence $\tilde{M}_{\tilde{\phi}}(\lambda(g)) = M_\phi( \lambda(g) )$, where $M_\phi$ is as defined in Theorem \ref{thm_main_theorem}. Applying the theorem one obtains that $\tilde{\phi}$ is a Herz-Schur-multiplier, and $\| \tilde{\phi} \|_{HS} \le \| \phi \|_\mathcal{C}$.
\end{proof}

\subsection{Relation to Houdayer and Ricard's results}
\label{radial_subsection_RH_examples}
Recently, C. Houdayer and {\'E}. Ricard in \cite{HoudayerRicard2011} proved results concerning radial multipliers on free Araki-Woods factors related to functions from a class $\mathcal{C}'$ quite similar to $\mathcal{C}$. Although their results apply to different objects than those considered in this paper, we will discuss in this section the issue of how their methods could be applied to prove Theorem \ref{thm_main_theorem}.

We start by defining the class of functions $\mathcal{C}'$, mentioned above.

\begin{definition}
\label{def_class_Cprime}
Let $\mathcal{C}'$ denote the set of functions $\phi: \mathbb{N}_0 \to \mathbb{C}$ for which the Hankel matrix $\hat{h} = (\phi(i+j) - \phi(i+j+2))_{i,j \ge 0}$ is of trace-class.
\end{definition}
Observe that this implies the existence of $c_1, c_2 \in \mathbb{C}$ and a unique $\psi: \mathbb{N}_0 \to \mathbb{C}$ such that $\phi(n) = c_1 + (-1)^n c_2 + \psi(n)$ and $\lim_{n \to \infty} \psi(n) = 0$.
For $\phi \in \mathcal{C}'$ put $\| \phi \|_{\mathcal{C}'} = |c_1| + |c_2| + \|\hat{h}\|_1$.

In \cite{HoudayerRicard2011} the following two results for functions in the class $\mathcal{C}'$ are proved, note the resemblance with Theorem \ref{thm_main_theorem}. In what follows $\mathcal{T}$ denotes the Toeplitz algebra, and furthermore, $\Gamma(H,U_t)''$ denotes the free Araki-Woods factor associated to a real Hilbert space $H$ and a one parameter group of orthogonal transformations $(U_t)$. (See \cite[Sections 2.5 and 3.1]{HoudayerRicard2011} for more precise definitions).

\begin{theorem}[{\cite[Proposition 3.3]{HoudayerRicard2011}}]
A function $\phi$ belongs to $\mathcal{C}'$ if and only if the operator $\gamma$ defined by $\gamma(S^i (S^*)^j) = \phi(i+j)$ extends to a bounded map on $\mathcal{T}$.
Moreover, $\| \gamma \|_{\mathcal{T}^*} = \| \phi \|_{\mathcal{C}'}$, and we say that $\gamma$ is the \emph{radial functional} associated with $\phi$.
\end{theorem}

\begin{theorem}[{\cite[Theorem 3.5]{HoudayerRicard2011}}]
Let $\phi: \mathbb{N}_0 \to \mathbb{C}$. Then $\phi$ defines a completely bounded radial multiplier on $\Gamma(H,U_t)''$ if and only if the radial functional $\gamma$ on $\mathcal{T}$ associated to $\phi$ is bounded. Moreover, $\| M_\phi \|_{cb} = \| \gamma \|_{\mathcal{T}^*}$.
\end{theorem}

A similar argument as in the proof of these theorems could be used to prove Theorem \ref{thm_main_theorem}, if one could prove the existence of a $*$-isomorphism
\begin{align}
\pi: C^*(L_\gamma | \gamma \in \Lambda(1)) &\to C^*(L_\gamma | \gamma \in \Lambda(1)) \otimes C^*(S^2,SS^*)
\end{align}
such that
\begin{align}
\label{eq_houdayer_ricard_pi}
\pi( L_\xi L_\eta^* ) &= \left\{ \begin{array}{ll}
L_\xi L_\eta^* \otimes S^{2k} (S^*)^{2l} & \qquad \text{in Case 1} \\
L_\xi L_\eta^* \otimes S^{2k-1} (S^*)^{2l-1} & \qquad \text{in Case 2}
\end{array} \right.
\end{align}
for all $k,l \ge 0$ and $\xi \in \Lambda(k), \eta \in \Lambda(l)$.

Indeed, if this were the case we could choose $w \in C^*(S^2,SS^*)^*$ as
\begin{align}
w(S^k (S^*)^l) = \left\{ \begin{array}{ll}
\phi \left(\frac{k+l}{2}\right) &\qquad \text{if } k+l \text{ even} \\
0 & \qquad \text{otherwise}.
\end{array} \right.
\end{align}

This functional would be bounded if $\hat{h}$ were trace-class. Moreover $\| Id \otimes w \|_{cb} = \| w \| = \| \phi \|_\mathcal{C}$.
Letting $T$ be defined by $T = (Id \otimes w) \circ \pi$ we would get
\begin{align}
T( L_\xi L_\eta^* ) = \left\{ \begin{array}{ll}
\phi(k+l) L_\xi L_\eta^* & \qquad \text{in Case 1} \\
\phi(k+l-1) L_\xi L_\eta^*  & \qquad \text{in Case 2}
\end{array} \right.
\end{align}
with $\| T \|_{cb} \le \| w \| \| \pi \|_{cb} \le \| w \| = \| \phi \|_\mathcal{C}$.
Hence, by Lemma \ref{lem_equiv_T_mult_T_cases} $T$ would be $M_\phi$ as defined in Theorem \ref{thm_main_theorem} and thus be completely bounded.

It is however not possible to construct such an isomorphism. Let for instance $|I|=1$ and $\dim(H)=2$ and let $e_{ij}$ denote the matrix units with respect to the basis $(\Omega,\gamma)$ of $H$. Then we have $e_{01} e_{10} = e_{00}$, but $\Phi(e_{01}) \Phi(e_{10}) = e_{00} \otimes S^*S \neq 1 \otimes 1 - e_{11} \otimes SS^* = \Phi(e_{00})$.

However, note that it would be sufficient if there existed a unital completely positive $\pi$ satisfying \eqref{eq_houdayer_ricard_pi}.
To find such an operator we can regard $l^2(\mathbb{N}_0) = l^2(\mathbb{N}_0)^{\text{even}} \oplus l^2(\mathbb{N}_0)^{\text{odd}}$. In this case $S^2$ on $l^2(\mathbb{N}_0)$ can be realized as $S \oplus S$ and $S S^*$ on $l^2(\mathbb{N}_0)$ can be realized as $S S^* \oplus 1$.
Then it would be enough to find unital completely positive operators $\pi_1, \pi_2$ such that
\begin{align}
\label{eq_example_RH_T1_requirements}
\pi_1( L_\xi L_\eta^* ) & = L_\xi L_\eta^* \otimes S^k (S^*)^l,
\end{align}
respectively,
\begin{align}
\label{eq_example_RH_T2_requirements}
\pi_2( L_\xi L_\eta^* ) & = \left\{ \begin{array}{ll}
L_\xi L_\eta^* \otimes S^{k} (S^*)^{l} & \qquad \text{in Case 1} \\
L_\xi L_\eta^* \otimes S^{k-1} (S^*)^{l-1} & \qquad \text{in Case 2}.
\end{array} \right.
\end{align}

Since $\pi(1)=\pi(2)=1$ we have $\|\pi_1\|_{cb},\|\pi_2\|_{cb} \le 1$ and $\pi = \pi_1 \oplus \pi_2$ is unital completely positive too. Hence $T = (Id \otimes w) \circ \pi$ would be as desired.

Set $U_n = \sum_{i=0}^\infty P_{i+n} \otimes e_{i0}$, where $e_{ij}$ are the matrix units in $B(l^2(\mathbb{N}_0))$, and use the convention $P_m = 0$ if $m<0$. Now it can be shown that
\begin{align}
\pi_1(x) & = \sum_{n=-\infty}^0 U_n (x \otimes 1) U_n^* + \sum_{n=1}^\infty U_n ( \rho^n(x) \otimes 1 ) U_n^*,
\end{align}
respectively,
\begin{align}
\pi_2(x) & = \sum_{n=-\infty}^0 U_n (x \otimes 1) U_n^* + \sum_{n=1}^\infty U_n ( \rho^{n-1}(\epsilon(x)) \otimes 1 ) U_n^*
\end{align}
are unital completely positive and fulfill \eqref{eq_example_RH_T1_requirements} and \eqref{eq_example_RH_T2_requirements}. The proof of this fact can be given by an argument quite similar to that given in Sections \ref{section_technical_lemmas} and \ref{section_main_proof}. We leave the details to the reader.

\section{\texorpdfstring{Integral representation of functions from $\mathcal{C}$}{Integral representation of functions from C}}
\label{section_integral}

In \cite{HaagerupSteenstrupSzwarc2010} the following integral representation was proved for $\phi \in \mathcal{C}'$, where $\mathcal{C}'$ is the set of functions $\phi: \mathbb{N}_0 \to \mathbb{C}$ from Definition \ref{def_class_Cprime}. The set $\mathcal{C}'$ is not defined in \cite{HaagerupSteenstrupSzwarc2010}, but the result follows form \cite[Theorem 2.12 and Theorem 4.2]{HaagerupSteenstrupSzwarc2010}.

\begin{theorem}
\label{thm_HSS_212_42}
Let $\psi: \mathbb{N}_0 \to \mathbb{C}$ be a function. Then the following are equivalent:
\begin{enumerate}
\item $\psi \in \mathcal{C}'$
\item There exists a complex Borel measure $\mu$ on $\mathbb{D}$ and constants $c_+, c_- \in \mathbb{C}$, such that
\begin{align}
\label{eq_psi_integral}
\psi(n) &= c_+ + (-1)^n c_- + \int_\mathbb{D} s^n d\mu(s) < \infty
\end{align}
and
\begin{align}
\int_\mathbb{D} \frac{|1-s^2|}{1-|s|^2} d|\mu|(s) < \infty.
\end{align}
\end{enumerate}
Moreover, for $\phi \in \mathcal{C}'$, the measure $\mu$ in \eqref{eq_psi_integral} can be chosen such that
\begin{align}
|c_+| + |c_-| + \int_\mathbb{D} \frac{|1-s^2|}{1-|s|^2} d|\mu|(s) \le \frac{8}{\pi} \|\psi\|_{\mathcal{C}'}.
\end{align}
\end{theorem}

We will prove next a similar characterization of functions in $\mathcal{C}$:

\begin{theorem}
\label{thm_phi_C_mu}
Let $\phi: \mathbb{N}_0 \to \mathbb{C}$ be a function. Then the following are equivalent:
\begin{enumerate}
\item $\phi \in \mathcal{C}$
\item There exists a constant $c \in \mathbb{C}$ and a complex Borel measure $\nu$ on $\mathbb{D}$ such that
\begin{align}
\label{eq_phi_integral}
\phi(n) &= c + \int_\mathbb{D} s^n d\nu(s)
\end{align}
and
\begin{align}
\int_\mathbb{D} \frac{|1-s|}{1-|s|} d|\nu|(s) < \infty.
\end{align}
\end{enumerate}
Moreover, for $\phi \in \mathcal{C}$, the measure $\nu$ in \eqref{eq_phi_integral} can be chosen such that
\begin{align}
|c| + \int_\mathbb{D} \frac{|1-s|}{1-|s|} d|\nu|(s) \le \frac{8}{\pi} \| \phi \|_\mathcal{C}.
\end{align}
\end{theorem}

\begin{proof}
(1) implies (2). Let $\phi \in \mathcal{C}$ and put
\begin{align}
\tilde{\phi}(n) = \left\{ \begin{array}{cl} \phi(\frac{n}{2}) & \text{ if $n$ is even} \\ 0 & \text{ if $n$ is odd.} \end{array} \right.
\end{align}
Then by Definition \ref{def_class_C} and Definition \ref{def_class_Cprime}, $\tilde{\phi} \in \mathcal{C}'$ and $\| \tilde{\phi} \|_{\mathcal{C}'} = \| \phi \|_\mathcal{C}$. From Theorem \ref{thm_HSS_212_42} there exists a complex measure $\mu$ on $\mathbb{D}$ and constants $c_+, c_- \in \mathbb{C}$ such that
\begin{align}
\psi(n) &= \tilde{\phi}(2n) = c_+ + c_- + \int_\mathbb{D} s^n d\mu(s) < \infty
\end{align}
and
\begin{align}
|c_+| + |c_-| + \int_\mathbb{D} \frac{|1-s^2|}{1-|s|^2} d|\mu|(s) \le \frac{8}{\pi} \|\psi\|_{\mathcal{C}}.
\end{align}

Let $\nu$ be the range measure of $\mu$ by the map $s \to s^2$ of $\mathbb{D}$ onto $\mathbb{D}$, and put $c = c_+ + c_-$. Then $|\nu|$ is less or equal to the range measure of $|\mu|$ by the map $s \to s^2$. Hence
\begin{align}
\phi(n) = c + \int_\mathbb{D} s^{2n} d\mu(s) = c + \int_\mathbb{D} s^n d\nu(s)
\end{align}
and
\begin{align}
|c| + \int_\mathbb{D} \frac{|1-s|}{1-|s|} d|\nu|(s) \le |c_+| + |c_-| + \int_\mathbb{D} \frac{|1-s^2|}{1-|s^2|} d|\mu|(s) \le \frac{8}{\pi} \|\phi\|_\mathcal{C}.
\end{align}
This proves (1) implies (2) and the last statement in Theorem \ref{thm_phi_C_mu}.

Conversely if (2) holds, the Hankel matrices $h, k$ from Definition \ref{def_class_C} have the entries
\begin{align}
h_{ij} &= \int_\mathbb{D} s^{i+j}(1-s) d\nu(s)
\end{align}
and
\begin{align}
k_{ij} &= \int_\mathbb{D} s^{i+j}s(1-s) d\nu(s).
\end{align}
By the proof of Corollary \ref{cor_phi_s},
\begin{align}
\|(s^{i+j})_{i,j \ge 0}\|_1 &= \frac{1}{1-|s|^2}, \qquad s \in \mathbb{D}.
\end{align}
Hence
\begin{align}
\|h\|_1 + \|k\|_1 \le \int_\mathbb{D} \frac{|1-s|+|s(1-s)|}{1-|s|^2} d|\nu|(s) = \int_\mathbb{D} \frac{|1-s|}{1-|s|} d|\nu|(s) < \infty
\end{align}
which shows that $\phi \in \mathcal{C}$.
\end{proof}

In \cite{Ozawa2008} N. Ozawa proved that if $\Gamma$ is a discrete hyperbolic group (in the sense of M. Gromov \cite{Gromov1987}), then $\Gamma$ is weakly amenable. The proof was obtained by showing that the metric $d: \Gamma \times \Gamma \to \mathbb{N}_0$ (w.r.t. the Cayley graph of $\Gamma$) satisfies three properties (1), (2) and (3) listed in \cite[Theorem 1]{Ozawa2008}. 

As an application of Theorem \ref{thm_phi_C_mu}, we will show below, that the first condition (1) from \cite{Ozawa2008} is sufficient to prove that $\Gamma$ is weakly amenable. For the definition of weak amenability and of the constant $\Lambda(\Gamma)$ for a weakly amenable group $\Gamma$, we refer to \cite[Section 12.3]{BrownOzawa2008}.

Recall that a metric on a discrete metric space $(X,d)$ is called \emph{proper} if the ball $B(x,r) = \{ y \in X: d(x,y) < r \}$ is finite for all $x \in X$ and all $r > 0$.

\begin{theorem}
\label{thm_sd_weakly_amenable}
Let $\Gamma$ be a discrete countable group and let $d: \Gamma \times \Gamma \to \mathbb{N}_0$ be a proper left invariant metric. Put
\begin{align}
\phi_s(x) = s^{d(x,e)}, \qquad s \in \mathbb{D}, x \in \Gamma.
\end{align}
Assume that there exists a constant $C \ge 1$, such that $\psi_s \in M_0A(\Gamma)$ for all $s \in \mathbb{D}$ and
\begin{align}
\| \phi_s \|_{M_0A(\Gamma)} \le C \frac{|1-s|}{1-|s|}, \qquad s \in \mathbb{D}.
\end{align}
Then $\Gamma$ is weakly amenable with constant $\Lambda(\Gamma) \le C$.
\end{theorem}

\begin{remark}
As in \cite{HaagerupSteenstrupSzwarc2010} we have used the notation $M_0A(\Gamma)$ for the set of completely bounded Fourier multipliers on $\Gamma$. Note that in \cite[Section 12.3]{BrownOzawa2008} the space $M_0A(\Gamma)$ is denoted $B_2(\Gamma)$.
\end{remark}

We first prove
\begin{lemma}
\label{lem_chi_n}$\quad$\newline
\begin{enumerate}
\item Put $\chi_n(k) = \delta_{kn}$ for $n,k \ge 0$. Then $\chi_n \in \mathcal{C}$ and
\begin{align}
\| \chi_n \|_\mathcal{C} \le \max\{1, 4n\}, \qquad n \ge 0.
\end{align}
\item For  $r \in (0,1)$ and $l \ge 0$, put
\begin{align}
\phi_r(k) &= r^k \\
\phi_{r,n}(k) &= \left\{ \begin{array}{cl} r^k & 0 \le k \le n \\ 0 & k > n. \end{array} \right.
\end{align}
then $\phi_r, \phi_{r,n} \in \mathcal{C}$,  $\| \phi_r \|_\mathcal{C} = 1$ and for fixed $r \in (0,1)$
\begin{align}
\lim_{n \to \infty} \| \phi_r - \phi_{r,n} \|_\mathcal{C} = 0.
\end{align}
\end{enumerate}
\end{lemma}

\begin{proof}
From Definition \ref{def_class_C} we have $\chi_n \in \mathcal{C}$, and $\| \chi_n \|_\mathcal{C} = \| H_n \|_1 + \| K_n \|_1$ where
\begin{align}
H_n(i,j) &= \chi_n(i+j) - \chi_n(i+j+1) \\
K_n(i,j) &= \chi_n(i+j+1) - \chi_n(i+j+2).
\end{align}

If $H=(h_{ij})_{i,j=0}^\infty$ is a matrix of complex numbers for which $\sum_{i,j} |h_{ij}| < \infty$, then $H$ is of trace class and $\| H \|_1 \le \sum_{i,j=0}^\infty | h_{ij} |$.
Hence $\| \chi_n \|_\mathcal{C} \le (2n+1) + (2n-1)$ for $n \ge 1$ and $\| \chi_0 \|_\mathcal{C} \le 1$ which proves (1). It follows from Corollary \ref{cor_phi_s}, that $\|\phi_r\|_\mathcal{C} = 1$, $0 < r < 1$. By (1),
\begin{align}
\| \phi_r - \phi_{r,n} \|_\mathcal{C} = \| \sum_{k=n+1}^\infty r^k \chi_k \|_\mathcal{C} \le \sum_{k = n+1}^\infty 4kr^k
\end{align}
which proves (2).
\end{proof}

\begin{proof}[Proof of Theorem \ref{thm_sd_weakly_amenable}]
Let $\phi: \mathbb{N}_0 \to \mathbb{C}$ be a function from $\mathcal{C}$, and put
\begin{align}
\tilde{\phi}(x) = \phi(d(x,e)), \qquad x \in \Gamma.
\end{align}
Then by \eqref{eq_psi_integral} and the integral representation of $\phi$ from Theorem \ref{thm_phi_C_mu} it follows that $\tilde{\phi}$ is a completely bounded Fourier multiplier on $\Gamma$ and that
\begin{align}
\label{eq_tilde_phi_M0A}
\| \tilde{\phi} \|_{M_0A(\Gamma)} \le \frac{8}{\pi} \| \phi \|_\mathcal{C}.
\end{align}

Let $\phi_r$ and $\phi_{r,n}$ be as in Lemma \ref{lem_chi_n}. Then by \eqref{eq_psi_integral} $\| \tilde{\phi}_r \|_{M_0A(\Gamma)} \le C$.

Moreover by \eqref{eq_tilde_phi_M0A} and Lemma \ref{lem_chi_n}
\begin{align}
\lim_{n \to \infty} \| \tilde{\phi}_r - \tilde{\phi}_{r,n} \|_{M_0A(\Gamma)} = 0
\end{align}
for fixed $r \in (0,1)$. Put $r_k = 1 - 1/k$, $k \ge 1$ and chose for each $k \ge 2$ an $n_k \ge k$, such that
\begin{align}
\| \tilde{\phi}_{r_k} - \tilde{\phi}_{r_k,n_k} \|_{M_0A(\Gamma)} \le \frac{1}{k}.
\end{align}

Then $\psi_k = \tilde{\phi}_{r_k,n_k}$ form a sequence of finitely supported functions on $\Gamma$, such that $\| \psi_k\|_{M_0A(\Gamma)} < C + 1/k$, and $\lim_{k \to \infty} \psi_k(x) = 1$ for all $x \in \Gamma$. Hence $\Gamma$ is weakly amenable and $\Lambda(\Gamma) \le C$.
\end{proof}

\begin{remark}
\label{remark_phi_tilde_Fourier}
By \cite[Theorem 1]{Ozawa2008} and the proof of Theorem \ref{thm_sd_weakly_amenable} it follows that for every hyperbolic group $\Gamma$ and every $\phi \in \mathcal{C}'$, the function
\begin{align}
\tilde{\phi}(x) = \phi(d(x,e)), \qquad x \in \Gamma
\end{align}
is a completely bounded Fourier multiplier on $\Gamma$ and $\| \tilde{\phi} \|_{M_0A(\Gamma)} \le \frac{8C}{\pi} \| \phi \|_{\mathcal{C}'}$, where $C$ is the constant in \cite[Theorem 1 (1)]{Ozawa2008}.
\end{remark}

\bibliography{qualifying}

\end{document}